\documentclass[11pt]{amsart}
\usepackage{amsmath}
\usepackage{amsxtra}
\usepackage{amscd}
\usepackage{amsthm}
\usepackage{amsfonts}
\usepackage{amssymb}
\usepackage{eucal}

\newtheorem{thm}{Theorem}[section]

\newtheorem{cor}[thm]{Corollary}
\newtheorem{lem}[thm]{Lemma}
\newtheorem{prop}[thm]{Proposition}

\theoremstyle{remark}
\newtheorem{remark}[thm]{Remark}

\theoremstyle{definition}

\numberwithin{equation}{section}

\newcommand{\Ref}[1]{{$($\ref{#1}$)$}}
\newcommand{\bean}{\begin{eqnarray}}
\newcommand{\eean}{\end{eqnarray}}
\newcommand{\be}{\begin{displaymath}}
\newcommand{\ee}{\end{displaymath}}
\newcommand{\bea}{\begin{eqnarray*}}   
\newcommand{\eea}{\end{eqnarray*}}

\newcommand{\thmref}[1]{Theorem~\ref{#1}}
\newcommand{\secref}[1]{Section~\ref{#1}}
\newcommand{\lemref}[1]{Lemma~\ref{#1}}
\newcommand{\propref}[1]{Proposition~\ref{#1}}
\newcommand{\corref}[1]{Corollary~\ref{#1}}
\newcommand{\remref}[1]{Remark~\ref{#1}}

\newcommand{\nc}{\newcommand}
\nc{\on}{\operatorname}
\nc{\ch}{\mbox{ch}}
\nc{\Z}{{\mathbb Z}}
\nc{\C}{{\mathbb C}}
\nc{\pone}{{\mathbb C}{\mathbb P}^1}
\nc{\pa}{\partial}
\nc{\F}{{\mathcal F}}
\nc{\arr}{\rightarrow}
\nc{\larr}{\longrightarrow}
\nc{\al}{\alpha}
\nc{\ri}{\rangle}
\nc{\lef}{\langle}
\nc{\W}{{\mathcal W}}
\nc{\la}{\lambda}
\nc{\ep}{\epsilon}
\nc{\su}{\widehat{{\mathfrak sl}}_2}
\nc{\sw}{{\mathfrak s}{\mathfrak l}}
\nc{\g}{{\mathfrak g}}
\nc{\h}{{\mathfrak h}}
\nc{\n}{{\mathfrak n}}
\nc{\N}{\widehat{\n}}
\nc{\G}{\widehat{\g}}
\nc{\De}{\Delta_+}
\nc{\gt}{\widetilde{\g}}
\nc{\Ga}{\Gamma}
\nc{\one}{{\mathbf 1}}
\nc{\z}{{\mathfrak Z}}
\nc{\zz}{{\mathcal Z}}
\nc{\Hh}{{\mathcal H}_\beta}
\nc{\qp}{q^{\frac{k}{2}}}
\nc{\qm}{q^{-\frac{k}{2}}}
\nc{\La}{\Lambda}
\nc{\wt}{\widetilde}
\nc{\wh}{\widehat}
\nc{\qn}{\frac{[m]_q^2}{[2m]_q}}
\nc{\cri}{_{\on{cr}}}
\nc{\kk}{h^\vee}
\nc{\sun}{\widehat{\sw}_N}
\nc{\hh}{{\mathbf H}_{q,t}}
\nc{\HH}{{\mathcal H}_{q,t}}
\nc{\hhh}{{\mathcal H}_{q,1}}
\nc{\ca}{\wt{{\mathcal A}}_{h,k}(\sw_2)}
\nc{\si}{\sigma}
\nc{\gl}{\widehat{{\mathfrak g}{\mathfrak l}}_2}
\nc{\el}{\ell}
\nc{\s}{T}
\nc{\bi}{\bibitem}
\nc{\om}{\omega}
\nc{\WW}{\W_\beta}
\nc{\scr}{{\mathbf S}}
\nc{\ab}{{\mathbf a}}
\nc{\rr}{r}
\nc{\ol}{\overline}
\nc{\con}{qt^{-1} + q^{-1}t}
\nc{\den}{q^{\el-1} t^{-\el+1}+ q^{-\el+1} t^{\el-1}}
\nc{\ds}{\displaystyle}
\nc{\B}{B}
\nc{\A}{A^{(2)}_{2\el}}
\nc{\GG}{{\mathcal G}}
\nc{\UU}{{\mathcal U}}
\nc{\MM}{{\mathcal M}}
\nc{\CC}{{\mathcal C}}
\nc{\GL}{^L\G}
\nc{\gL}{^L\g}
\nc{\dzz}{\frac{dz}{z}}
\nc{\Res}{\on{Res}}
\nc{\rep}{{\mathcal R}ep \;}
\nc{\uqg}{U_q \G}
\nc{\ueg}{U^{\on{res}}_\epsilon \G}
\nc{\tueg}{\wt{U}_\epsilon \G}
\nc{\fueg}{U^{\on{fin}}_\epsilon \G}
\nc{\uegg}{U^{\on{res}}_\epsilon \g}
\nc{\dueg}{\dot{U}^{\on{res}}_\epsilon \G}
\nc{\uqgr}{U_q^{\rm res}\G}
\nc{\uqgg}{U_q \g}
\nc{\uqggr}{U_q^{\rm res}\g}
\nc{\duqgr}{\dot{U}^{\on{res}}_q \G}
\nc{\duqg}{\dot{U}_q \G}
\nc{\duegs}{\dot{U}^{\on{res}}_{\epsilon^*} \G}
\nc{\duqges}{{}^*\dot{U}_\ep \G}
\nc{\uegs}{U^{\on{res}}_{\epsilon^*} \G}
\nc{\mc}{\mathcal}
\nc{\Cal}{\mathcal}
\nc{\calp}{{\Cal P}}
\nc{\bp}{{\mathbf P}}
\nc{\bq}{{\mathbf Q}}
\nc{\bb}{{\mathfrak b}}
\nc{\uqb}{U_q \bb_-}
\nc{\uqn}{U_q \wt{{\mathfrak n}}}
\nc{\uqh}{U_q \wt{{\mathfrak h}}}
\nc{\uqhh}{U_q \wh{{\mathfrak h}}}
\nc{\uqnn}{U_q \wh{{\mathfrak n}}}
\nc{\ot}{\otimes}
\nc{\R}{{\mc R}}
\nc{\uqbb}{U_q \wt{\g}}
\nc{\yy}{{\mc Y}}
\nc{\uqsl}{U_q \widehat{\sw}_2}
\nc{\ga}{\gamma}
\nc{\Ab}{{\mathbf A}}
\nc{\Yb}{{\mathbf Y}}
\nc{\yb}{{\mathbf y}}
\nc{\uh}{U \wt{{\mathfrak h}}}
\nc{\uhh}{U \wh{{\mathfrak h}}}

\nc{\us}{\underset}
\nc{\opl}{\oplus}
\nc{\yyy}{\wh{\yy}}
\nc{\ovl}{\overline}
\nc{\beq}{\begin{equation}}
\nc{\Fq}{{\mathcal F}}
\nc{\Mq}{{\mathcal M}}
\nc{\Rep}{\on{Rep}}
\nc{\sssec}{\subsubsection}
\nc{\ssec}{\subsection}
\nc{\lan}{\langle}
\nc{\ran}{\rangle}

\nc{\uqhJ}{U_q \widehat{\h}_J^\perp}
\nc{\uqsli}{U_{q_i} \widehat{\sw}_2}
\nc{\len}{\on{length}}

\def\bin[#1;#2]{\left[\begin{matrix}{\displaystyle
#1}\\{\displaystyle #2}\end{matrix}\right]}

\def\qbin2[#1;#2;#3]{\left[\begin{matrix}{\displaystyle
#1;\displaystyle #2}\\{\displaystyle #3}\end{matrix}\right]}
 
\nc{\ghat}{\widehat{\g}}
\nc{\umod}{\dot{\mathbf U}}

\begin{document}

\title[The $q$--characters at roots of unity]{The $q$--characters at
roots of unity}

\author {Edward Frenkel}

\author{Evgeny Mukhin}

\address{Department of Mathematics, University of California, Berkeley, CA
94720, USA}

\begin{abstract}
We consider various specializations of the untwisted quantum affine
algebras at roots of unity. We define and study the $q$--characters of
their finite-dimensional representations.
\end{abstract}

\maketitle

\section{Introduction}

The theory of $q$--characters of finite-dimensional representations of
quantum affine algebras was developed in \cite{FR,FM}. In those works,
$q$ was assumed to be a generic non-zero complex number (i.e., not a
root of unity). In the present paper we extend the results of
\cite{FR,FM} to the case when $q$ is a root of unity.

There are different versions of the quantum affine algebra $\uqg$ when
$q$ is specialized to a root of unity $\ep$: the non-restricted
specialization $\tueg$ studied by Beck and Kac \cite{BK}, the
restricted specialization $\ueg$ studied by Chari and Pressley
\cite{CP:root} (following the general definition due to Lusztig
\cite{L}), and the ``small'' affine quantum group
$U^{\on{fin}}_\ep \G$, which is the image of the natural homomorphism
$\tueg \to \ueg$.

It was shown in \cite{CP:root} that all irreducible finite-dimensional
representations of $\ueg$ (of type 1) are highest weight
representations with respect to a triangular decomposition in terms of
the Drinfeld generators (more precisely, these results were obtained
in \cite{CP:root} in the case when $\ep$ is a root of unity of odd
order, but we extend them here to all roots of unity). Using these
results, we define the $\ep$--character of a $\ueg$--module $V$ via
the generalized eigenvalues of a commutative subalgebra of $\ueg$ on
$V$. We establish various properties of the $\ep$--character
homomorphism; in particular, we show that the $q$--characters
specialize to $\ep$--characters as $q \to \ep$.

In addition, using results of \cite{BK} we show that all
irreducible finite-dimensional representations of $\fueg$ (of type 1)
are highest weight representations, describe their highest weights and
define the corresponding $\ep$--characters.  We show that each
irreducible $\fueg$--module admits a unique structure of
$\ueg$--module. This allows us to obtain the $\ep$--characters of
$\fueg$--modules from the $\ep$--characters of $\ueg$--modules.

Finally, we study the quantum Frobenius homomorphism, following
Lusztig's definition \cite{L}. Let $l$ be the order of $\ep^2$ and
$\ep^* =\ep^{l^2}\in\{ \pm 1\}$. The Frobenius homomorphism maps
$\ueg$ (more precisely, its modified version) either to
$U_{\ep^*}^{\on{res}} \G$ or to $U_{\ep^*}^{\on{res}} \G^L$,
depending on whether $l$ is divisible by the lacing number $r^\vee$ or
not. Here $\G^L$ denotes the Langlands dual Lie algebra to $\G$ (it is
twisted if $\G$ is not simply-laced) and we use the notation $\Rep U$
for the Grothendieck ring of finite-dimensional representations of an
algebra $U$. In this paper we deal only with the untwisted quantum
affine algebras, and so we consider the Frobenius homomorphism in
detail only in the case when $(l,r^\vee)=1$.

The Frobenius pull-backs of irreducible finite-dimensional
$U_{\ep^*}^{\on{res}} \G$--modules are irreducible $\ueg$--modules. We
describe their $\ep$--characters in terms of the ordinary
characters of irreducible finite-dimensional representations of the
simple Lie algebra $\g$.

A decomposition theorem (previously proved by Chari and Pressley
\cite{CP:root} in the case when $\ep$ has odd order) allows one to
decompose any irreducible finite-dimensional representation of $\ueg$
as a tensor product of a Frobenius pull-back and a module which
remains irreducible when restricted to $U^{\on{fin}}_\ep \G \subset
\ueg$. Thus, we obtain an isomorphism of vector spaces
$$
\Rep \ueg \simeq \Rep \fueg \otimes \Rep U_{\ep^*}^{\on{res}} \G.
$$

Similar results can be obtained for the twisted quantum affine
algebras. We will describe the $q$--characters and
$\ep$--characters of twisted quantum affine algebras in a separate
publication.
\medskip

\noindent{\bf Acknowledgments.} We thank V.~Chari, D.~Gaitsgory,
A.~Kirillov, Jr., and N.~Resheti\-khin for useful discussions. This
research was supported in part by the Packard Foundation and the NSF.

\section{Definitions}

\subsection{Root data}    \label{cartan}

Let $\g$ be a simple Lie algebra of rank $\el$. Let $( \cdot,\cdot )$
be the invariant inner product on $\g$, such that the square of the
length of the maximal coroot $\al^\vee_{\on{max}}$ equals $2$. The
induced inner product on the dual space to the Cartan subalgebra $\h$
of $\g$ is also denoted by $( \cdot,\cdot )$. Denote by $I$ the set
$\{ 1,\ldots,\el \}$. Let $\{ \al_i \}_{i \in I }$, $\{ \al^\vee_i
\}_{i \in I }$ $\{ \om_i \}_{i \in I}$ be the sets of simple roots,
simple coroots and fundamental weights of $\g$, respectively.  We
denote by $\al_{\on{max}}$ the maximal (positive) root of $\g$.


Let $r^\vee$ be the maximal number of edges connecting two vertices of
the Dynkin diagram of $\g$. Thus, $r^\vee=1$ for simply-laced $\g$,
$r^\vee=2$ for $B_\el, C_\el, F_4$, and $r^\vee=3$ for $G_2$.

Set $\rr_i = \frac{(\al_i,\al_i)}{2}$. All $\rr_i$'s are equal to $1$
for simply-laced $\g$, and $r_i$'s are equal to $1$ or
$r^\vee$ depending on whether
$\al_i$ is short or long, for non-simply laced $\g$.

Let $C = (C_{ij})_{i,j\in I}$ be the Cartan matrix of $\g$, $C_{ij} =
\frac{2(\al_i,\al_j)}{(\al_i,\al_i)}$.

Let $\hat I=\{0,1,\dots,\ell\}$ and
let $(C_{ij})_{i,j\in\hat I}$ be the Cartan matrix of $\G$. We set
$r_0 = r^\vee$.

We also fix a choice of a function $o:\;I\to\{\pm 1\}$
such that $o(i)\neq o(j)$ whenever $(\al_i,\al_j) \neq 0$. 

\subsection{Quantum affine algebras}    \label{qaal}

Let $q$ be an indeterminate, $\C(q)$ the field of rational functions
of $q$ and $\C[q,q^{-1}]$ the ring of Laurent polynomials with complex
coefficients. For $n\in\Z$, $m\in\Z_{\geq 0}$ set \be
[n]_q=\frac{q^n-q^{-n}}{q-q^{-1}},\qquad \left[\begin{matrix}n
\\m\end{matrix}\right]_q=\frac{\prod_{i=n-m+1}^n [i]_q}{
\prod_{i=1}^m[i]_q}, \qquad [m]_q!=\prod_{i=1}^m[i]_q.
  \ee Let $q_i=q^{r_i}, i \in \hat{I}$.

In this paper we deal exclusively with finite-dimensional
representations of quantum affine algebras, all of which have level
zero. The (multiplicative) central element of a quantum affine
algebra acts as the identity on such representations. To simplify our
formulas, we will impose in the definition of the quantum affine
algebra the additional relation that this central element is equal to
$1$.

The quantum affine algebra $\uqg$ (of level zero) in the
Drinfeld-Jimbo realization \cite{Dr1,J} is an associative algebra over
$\C(q)$ with generators $x_i^{{}\pm{}}$ ($i\in\hat I$), $k_i^{{}\pm
1}$ ($i \in I$), and relations:
\begin{align*}
k_ik_i^{-1} = k_i^{-1}k_i &=1,\quad \quad k_ik_j =k_jk_i,\\
k_ix_j^{{}\pm{}}k_i^{-1} &= q_i^{{}\pm C_{ij}}x_j^{{}\pm},\\ [x_i^+ ,
x_j^-] &= \delta_{ij}\frac{k_i - k_i^{-1}}{q_i -q_i^{-1}},\\
\sum_{r=0}^{1-C_{ij}}(-1)^r \left[\begin{array}{cc} 1-C_{ij} \\ r \end{array}
\right]_{q_i} &
(x_i^{{}\pm{}})^rx_j^{{}\pm{}}(x_i^{{}\pm{}})^{1-C_{ij}-r} =0 \ \ \
\ (i\ne j).
\end{align*}
Here we use the notation $k_0 = \prod_{i\in I} k_i^{-a_i}$, where $\al_{\max} =
\sum_{i \in I} a_i \al_i$.

The algebra $\uqg$ has a structure of a Hopf algebra with the
comultiplication $\Delta$ and the antipode $S$ given on the generators
by the formulas: \bea \Delta(k_i) &=& k_i \otimes k_i, \\
\Delta(x^+_i) &=& x^+_i \otimes 1 + k_i \otimes x^+_i,\\ \Delta(x^-_i)
&=& x^-_i \otimes k_i^{-1} + 1 \otimes x^-_i, \eea \be S(x_i^+) =
-x_i^+ k_i,\qquad S(x_i^-) = - k_i^{-1} x_i^-,\qquad S(k_i^{\pm 1}) =
k_i^{\mp 1}.  \ee



The Hopf algebra $\uqgg$ is defined as the $\C(q)$--subalgebra of $\uqg$ with
generators $x_i^{{}\pm{}}$, $k_i^{{}\pm 1}$, where $i\in I$.


Let $\sigma$ be an automorphism of the affine Dynkin diagram, i.e.,
$\sigma:\; \hat I\to \hat I$, such that $(\al_i,\al_j)=(
\sigma\al_i,\sigma\al_j)$. Then we also denote by $\sigma$ an
$\uqg$ automorphism defined by the formulas \bean\label{diag aut}
\sigma(x^{\pm}_i)=x^{\pm}_{\sigma(i)},\qquad
\sigma(k_i)=k_{\sigma(i)}\qquad (i\in\hat I).  \eean

Let $T_i$, $i\in \hat I$, be the $\C(q)$--algebra automorphisms of
$\uqg$, defined by the formulas \bea T_i(k_j) =
k_i^{-a_{ij}}k_j,\qquad\qquad\qquad\qquad\qquad\qquad\\
T_i((x_i^+)^{(n)}) = (-1)^nq_i^{-n(n-1)}(x_i^-)^{(n)}k_i^n,\quad
T_i((x_i^-)^{(n)}) = (-1)^nq_i^{n(n-1)}k_i^{-n}(x_i^+)^{(n)},\\
T_i((x_j^+)^{(n)}) = \sum_{r=0}^{-na_{ij}}
(-1)^{r}q_i^{-r}(x_i^+)^{(-na_{ij}-r)}(x_j^+)^{(n)}(x_i^+)^{(r)}
\qquad (i\neq j),\\ T_i((x_j^-)^{(n)}) = \sum_{r=0}^{-na_{ij}}
(-1)^{r}q_i^{r}(x_i^-)^{(r)}(x_j^-)^{(n)}(x_i^-)^{(-na_{ij}-r)} \qquad
(i\neq j).  \eea Here we use the following notation for divided powers
\be (x^{\pm}_i)^{(r)}=\frac{(x_i^\pm)^r}{[r]_{q_i}!}.  \ee If $i\in I$
then $T_i$ induces an $\C(q)$--automorphism of subalgebra $\uqgg$ of
$\uqg$, also denoted by $T_i$.  The automorphisms $T_i$ are the
operators $T''_{i,-1}$ introduced in \cite{L}, \S 41.1.2.


Next, we describe the Drinfeld ``new'' realization of $\uqg$.

\begin{thm}[\cite{Dr,KhT2,LSS,B}]    \label{defining}
The algebra $\uqg$ is isomorphic over $\C(q)$ to the algebra with
generators $x_{i,n}^{\pm}$ ($i\in I$, $n\in\Z$), $k_i^{\pm 1}$ ($i\in
I$), $h_{i,n}$ ($i\in I$, $n\in \Z\backslash 0$), with the following
relations:
\begin{align*}
k_ik_j = k_jk_i,\quad & k_ih_{j,n} =h_{j,n}k_i,\\
k_ix^\pm_{j,n}k_i^{-1} &= q_i^{\pm C_{ij}}x_{j,n}^{\pm},\\ [h_{i,n} ,
x_{j,m}^{\pm}] &= \pm \frac{1}{n} [n C_{ij}]_{q_i} x_{j,n+m}^{\pm},\\
x_{i,n+1}^{\pm}x_{j,m}^{\pm} -q_i^{\pm
C_{ij}}x_{j,m}^{\pm}x_{i,n+1}^{\pm} &=q_i^{\pm
C_{ij}}x_{i,n}^{\pm}x_{j,m+1}^{\pm} -x_{j,m+1}^{\pm}x_{i,n}^{\pm},\\
[h_{i,n},h_{j,m}] &= 0,\\ [x_{i,n}^+ , x_{j,m}^-]=\delta_{ij} & \frac{
\phi_{i,n+m}^+ - \phi_{i,n+m}^-}{q_i - q_i^{-1}},\\
\sum_{\pi\in\Sigma_s}\sum_{k=0}^s(-1)^k\left[\begin{array}{cc} s \\ k
\end{array} \right]_{q_i} x_{i, n_{\pi(1)}}^{\pm}\ldots
x_{i,n_{\pi(k)}}^{\pm} & x_{j,m}^{\pm} x_{i, n_{\pi(k+1)}}^{\pm}\ldots
x_{i,n_{\pi(s)}}^{\pm} =0,\ \ s=1-C_{ij},
\end{align*}
for all sequences of integers $n_1,\ldots,n_s$, and $i\ne j$, where
$\Sigma_s$ is the symmetric group on $s$ letters, and
$\phi_{i,n}^{\pm}$'s are determined by the formula
\begin{equation}    \label{series}
\Phi_i^\pm(u) := \sum_{n=0}^{\infty}\phi_{i,\pm n}^{\pm}u^{\pm n} =
k_i^{\pm 1} \exp\left(\pm(q_i-q_i^{-1})\sum_{m=1}^{\infty}h_{i,\pm m}
u^{\pm m}\right).
\end{equation}
\end{thm}

Let $\uqg^\pm$ be the subalgebras of $\uqg$ generated by
$x_{i,n}^\pm$ ($i\in I$, $n\in\Z$). Let $\uqg^0$ be the subalgebra 
of $\uqg$ generated by $k_i^{\pm 1}$, $h_{i,n}^\pm$ ($i\in I$,
$n\in\Z\setminus 0$).

In fact, according to \cite{B} (see also \cite{CP:root,BCP}), the
isomorphism described in Theorem \ref{defining} has the form \be
x_{i,r}^\pm=o(i)^rT_{\check\om_i}^{\mp r}(x_i^{\pm}), \ee where
$T_{\check{\om_i}}$ is an automorphism of $\uqg$ given by a certain
composition of braid group automorphisms $T_i$ and a diagram
automorphism $\sigma$ of type \Ref{diag aut}. We remark that in
\cite{B,CP:root} a different normalization for the braid group action
is used. In this paper we follow the conventions of \cite{BCP}, \cite{L}.

For any $a \in \C^\times$, there is a Hopf algebra automorphism
$\tau_a$ of $\uqg$ defined on the generators by the following
formulas:
\bean\label{tau}
\tau_a(x_{i,n}^{\pm})=a^nx_{i,n}^{\pm}, \quad \quad
\tau_a(\phi_{i,n}^{\pm})=a^n\phi_{i,n}^{\pm},
\eean
for all $i \in I$, $n \in \Z$. Given a $\uqg$--module $V$ and $a \in
\C^\times$, we denote by $V(a)$ the pull-back of $V$ under $\tau_a$.

\subsection{Restricted integral form}

Let $\uqgr$ be the $\C[q,q^{-1}]$--subalgebra of $\uqg$ generated by
$k_i^{\pm1}$ and $(x_i^{\pm})^{(n)}$, $i\in \hat I$,
$n\in\Z_{>0}$. Then $\uqgr$ is a $\C[q,q^{-1}]$ Hopf subalgebra of
$\uqg$ preserved by automorphisms $T_i$, $i\in\hat I$.

Similarly, let $\uqggr$ be the $\C[q,q^{-1}]$--subalgebra of $\uqgg$
generated by 
$k_i^{\pm1}$ and $(x_i^{\pm})^{(n)}$, $i\in I$,
$n\in\Z_{>0}$. Then $\uqggr$ is a $\C[q,q^{-1}]$ Hopf subalgebra of
$\uqgg$ preserved by automorphisms $T_i$, $i\in I$.

For $i\in I$, $r\in\Z_{>0}$, define \be
\bin[k_i;r]=\prod_{s=1}^r\frac{k_iq_i^{1-s}-k_i^{-1}q_i^{s-1}}
{q_i^s-q_i^{-s}}.  \ee 
Note that $\bin[k_i;r]$ is denoted by $\qbin2[k_i;0;r]$ in \cite{CP:root}.

For $i\in I$, $n\in\Z$, define the elements $P_{i,n}\in \uqg$ by \be
{\mathcal P}_i^\pm(u)=\sum_{n=0}^\infty P_{i,\pm n} u^{\pm n}=
\exp\left(\mp\sum_{m=1}^\infty\frac{h_{i,\pm m}}{[m]_{q_i}}u^{\pm
m})\right).  \ee Note that \be\label{phi p} \Phi_i^\pm(u)=k_i^{\pm
1}\frac{{\mathcal P}_i^\pm(uq_i^{-1})}{{\mathcal P}_i^\pm(uq_i)}.  \ee

Denote by $\uqgr^\pm$ the $\C[q,q^{-1}]$--subalgebra of $\uqg$
generated by $(x_{i,r}^\pm)^{(n)}$, $i\in I$, $r\in\Z$,
$n\in\Z_{>0}$. Denote $\uqgr^0$ the $\C[q,q^{-1}]$--subalgebra of
$\uqg$ generated by $k_i$, $\bin[k_i;r]$, $i \in \hat{I}$ and
$P_{i,n}$, $i\in I$, $n\in\Z$, $r\in\Z_{>0}$.
 
\begin{thm}[\cite{CP:root}, Proposition 6.1]    \label{trdec}
We have $\uqgr^\pm\subset\uqgr$, $\uqgr^0\subset\uqgr$. 
The algebra $\uqgr$ is generated by the subalgebras $\uqgr^+$,
$\uqgr^-$ and $k_i^{\pm 1}$, $i\in I$.  
Moreover, the multiplication gives an isomorphism of vector spaces
\bean\label{triangular}
\uqgr\simeq\uqgr^-\otimes\uqgr^0\otimes\uqgr^+.  \eean
\end{thm}

Let $\epsilon\in\C$ be a primitive $s$-th root of unity, i.e.,
$\ep^s=1$ and $\ep^k\neq 1$ for $k=1,2,\dots,s-1$. 
For $i\in\hat{I}$, let $\ep_i=\ep^{r_i}$ and let
\be
l=\left\{\begin{matrix} s, & s \;\;{\rm is \;\;odd}\\ s/2, & s\;\;{\rm
is \;\;even}
\end{matrix}\right. \qquad
l_i=\left\{\begin{matrix} 
l, & l \;\;{\rm is \;\;not\;\; divisible\;\; by} \;\;r_i\\
s/r_i, &l \;\;{\rm is \;\; divisible \;\; by}\;\; r_i.
\end{matrix}\right. 
\ee 
In other words, $l$ is the order of $\ep^2$, and $l_i$ is the
order of $\ep_i^2$. Thus, $l$ is the minimal natural number with the
property $[l]_\ep=0$ and $l_i$ is the minimal natural number with the
property $[l_i]_{\ep_i}=0$.

Define the structure of a module over the ring $\C[q,q^{-1}]$ on $\C$
by the formula \be p(q,q^{-1}) \mapsto p(\epsilon,\ep^{-1}),\qquad
p(q,q^{-1})\in \C[q,q^{-1}]. \ee Denote this module by $\C_\ep$. Set
\bea \ueg=\uqgr\otimes_{\C[q,q^{-1}]}\C_\ep,\\
\uegg=\uqggr\otimes_{\C[q,q^{-1}]}\C_\ep.  \eea Also set \bea
\ueg^\pm=\uqgr^\pm\otimes_{\C[q,q^{-1}]}\C_\ep\subset \ueg,\\
\ueg^0=\uqgr^0\otimes_{\C[q,q^{-1}]}\C_\ep\subset \ueg.  \eea 
The algebra $\ueg$ is called the {\em restricted specialization} of
$\uqg$ at $q=\ep$.

\subsection{Non-restricted integral form and small affine quantum
group}

Let $\wt{U}_q \G$ be the $\C[q,q^{-1}]$--subalgebra of $\uqg$
generated by $k_i^{\pm1}$ and $x_i^{\pm}$, $i\in \hat I$.  This is a
$\C[q,q^{-1}]$ Hopf subalgebra of $\uqg$ preserved by automorphisms
$T_i$, $i\in\hat I$. Using \thmref{defining} we obtain that $\wt{U}_q
\G$ may also be described as the $\C[q,q^{-1}]$--subalgebra of $\uqg$
generated by $k_i^{\pm1}$, $x_{i,n}^{\pm}$ ($i\in I$,
$n\in\Z$), $h_{i,n}$ ($i\in I$, $n\in \Z\backslash 0$). Note that
$\wt{U}_q \G$ is a $\C[q,q^{-1}]$--subalgebra of $\uqgr$.

Define the {\em non-restricted specialization} of $\uqg$ at $q=\ep$ by
$$\tueg: = \wt{U}_q \G \otimes_{\C[q,q^{-1}]}\C_\ep.$$ We have a natural
homomorphism $\tueg \to \ueg$. Denote by $\fueg$ the image of this
homomorphism. We call $\fueg$ the {\em small affine quantum group} of
$\G$ at $q=\ep$. This is the subalgebra of $\ueg$ generated by
$k_i^{\pm1}$ and $x_i^{\pm}$, $i\in
\hat I$.

We remark that Chari and Pressley use in \cite{CP:root} the notation
$\fueg$ for a larger subalgebra of $\ueg$, which includes $P_{i,n}, i
\in I, n \in \Z$ in addition to the above generators.

Note that if $l=1$, we have: $\tueg=\fueg=\ueg$. Hence from now on,
when dealing with $\tueg$ or $\fueg$, we will always assume that $l>1$.

\begin{lem}    \label{rel in fin}
For $l>1$, the following relations hold in $\fueg$:
$(x_i^\pm)^{l_i}=0$ $(i \in \hat{I})$, $(x_{i,n}^\pm)^{l_i}=0$, and
$h_{i,ml_i}=0$ $(i \in I$, $n \in \Z$, $m \in \Z\backslash 0)$.
\end{lem}

\begin{proof}
This follows from the fact that the elements $(x_i^\pm)^{(l_i)},
(x_{i,n}^\pm)^{(l_i)}$, and $h_{i,ml_i}/[ml_i]_{q_i}$ belong to
$\ueg$.
\end{proof}

In particular, $\fueg$ is generated by $k_i^{\pm1}$,
$x_{i,n}^{\pm}$ ($i\in
I$, $n\in\Z$), $h_{i,n}$ ($i\in I$, $n\in \Z\backslash \Z l_i$).

\subsection{Finite-dimensional representations of $\ueg$}
\label{type 1}

A finite-dimensional representation of $\uqgg$ is said to be of type
$1$ if $k_i$, $i\in I$, act by semi-simple operators with eigenvalues
in $q_i^{\Z}$. A finite-dimensional representation $V$ of $\uegg$ is
said to be of type $1$ if $V$ has a basis of eigenvectors of the
elements $k_i$, $\bin[k_i;l_i], i\in I$, and for each
vector $v$ from this basis there exist integers $\la_i$, such that
$$
k_i v = \ep_i^{\la_i} v, \qquad \bin[k_i;l_i] v =
{\left[\begin{matrix}{\la_i}\\{ l_i}\end{matrix}\right]}_{\ep_i}v,
$$
(see \cite{CP}, \S 11.2). Such a vector $v$ is
called a weight vector of weight $\sum_{i\in I} \la_i \omega_i$.

A finite-dimensional $\uqg$--module (resp., $\ueg$--module) $V$ is
said to be of type $1$ if the restriction of $V$ to $\uqgg$ (resp.,
$\uegg$) is of type $1$.  In what follows we
consider only modules of type $1$.

It follows from the formulas defining the coproduct in $\uqg$ that the
tensor product of two type 1 representations is of type 1. Hence the
Grothendieck rings of the type 1 finite-dimensional representations of
$\uqg$ and $\ueg$ carry ring structures. We denote these rings by
$\on{Rep}\uqg$ and $\on{Rep}\ueg$, respectively.

A vector $v$ in a $\uqg$--module (resp., $\ueg$--module) $V$ is
called a highest weight vector if $\uqg^+v=0$ and $\uqg^0\cdot
v=\C(q)\cdot v$ (resp., $\ueg^+v=0$ and $\ueg^0\cdot
v=\C\cdot v$). If in addition $V=\uqg\cdot v$ (resp.,
$V=\ueg\cdot v$), then $V$ is called a highest weight representation.

\begin{thm}    \label{classe}
Every irreducible finite-dimensional representation of $\ueg$ is a
highest weight representation. An irreducible highest weight
representation $V$ with highest weight vector $v$ is
finite-dimensional if and only if there exists an $I$--tuple of
polynomials with constant term $1$, $${\bf P}=(P_i(u))_{i\in I},
\qquad P_i(u)=\prod_{j=1}^{j_i}(1-a_{ij}u), \qquad a_{ij}\in\C,$$ such
that \bea k_iv=\ep_i^{\deg P_i}v, \qquad
\bin[k_i;l_i]v={\left[\begin{matrix}{\deg P_i}\\{
l_i}\end{matrix}\right]}_{\ep_i}v,\\ \mc P^+_i(u)v=P_i(u)v,\qquad \mc
P^-_i(u)v=\wt{P}_i(u)v, \eea where
$\wt{P_i}(u)=\prod_{j=1}^{j_i}(1-(a_{ij}u)^{-1})$.
\end{thm}

\begin{proof}
In \cite{CP:root}, Proposition 8.1 and Theorem 8.2, the theorem is
proved when $s$ is odd and $(s,r^{\vee})=1$, using the triangular
decomposition \Ref{triangular} for $\uqgr$. The same proof applies in
general as the decomposition
\Ref{triangular} does not depend on the choice of $\ep$.
\end{proof}

The $I$--tuple of polynomials ${\bf P}$ is called the highest weight
of $V$. We write $V=V({\bf
P})$ for the irreducible representation $V$ with highest weight ${\bf
P}$. We refer to the polynomials $P_i(u)$ as Drinfeld polynomials.

There is an analogous description of irreducible finite-dimensional
$\uqg$--modules given as highest weight representations given in
\cite{CP,CP4}. Namely, irreducible $\uqg$--modules are also classified by
highest weights ${\bf P}$, which are $I$--tuples of polynomials in
$\C(q)[u]$ with constant coefficient $1$. We denote the irreducible
$\uqg$--module with highest weight ${\bf P}$ by $V({\bf P})_q$.

Let $\Lambda = \on{span}_{\Z} \{ \omega_1,\ldots,\omega_\ell \}$ be
the weight lattice of $\g$. Recall that $\Lambda$ carries a partial
order. Namely, for $\la,\la' \in \Lambda$, we say that $\la \geq \la'$
if $\la = \la' + \sum_{i \in I} n_i \al_i$, where $n_i\in \Z_{\geq
0}$. This induces a partial order on the set of $I$--tuples of
polynomials: we say that ${\bf P} \geq
{\bf Q}$ if $\sum_{i \in I} \omega_i\deg P_i \geq \sum_{i \in I} \omega_i\deg
Q_i$.

The next proposition shows that all irreducible finite-dimensional
representations of $\ueg$ can be obtained as subquotients of
specializations of $\uqg$--modules, and moreover, the corresponding
decomposition matrix is triangular with respect to the above partial
order.

\begin{prop}\label{spec rep}
Let $V({\bf P})_q$ be an irreducible highest weight $\uqg$--module of
dimension $d$ with highest vector $v$ and highest weight ${\bf P}$,
such that $P_i(u) \in \C[q,q^{-1},u]$ ($i \in I$). Then $V({\bf
P})_q^{\on{res}}:=\uqgr\cdot v$ is a free $\C[q,q^{-1}]$--module of
rank $d$, and $$V({\bf P})^{\on {res}}_\ep := V({\bf P})_q^{\on{res}}
\otimes_{\C[q,q^{-1}]} \C_\ep$$ is a $\ueg$--module of dimension
$d$. Moreover, the map \bea e:\;\on{Rep}\uqg&\to&\on{Rep}\ueg,\\
{}[V({\bf P})_q]{} &\mapsto& {}[V({\bf P})^{\on{res}}_\ep]{} \eea is a
surjective ring homomorphism, and
\begin{equation}    \label{multi}
[V({\bf P})^{\on{res}}_\ep] = [V({\bf P}_\ep)] + \sum_{{\bf Q}<{\bf
P}_\ep} m_{{\bf Q},{\bf P}_{\ep}} [V({\bf Q})], \qquad m_{{\bf Q},{\bf
P}_{\ep}} \in \Z_{\geq 0},
\end{equation}
where ${\bf P}_\ep$ is obtained from ${\bf P}$ by substituting
$q=\ep$.
\end{prop}

\begin{proof}
The module $V({\bf P})_q^{\on{res}}$ is a free $\C[q,q^{-1}]$--module
of rank $d$ by Lemma 4.6 i) in \cite{CP:weyl}.

Since $\uqgr$ is a Hopf subalgebra of $\uqg$, the map $e$ is a
homomorphism of rings. Since the degrees of the polynomials $P_i(u)$
determine the highest weight of $V({\bf P})_q$, considered as a
$\uqgg$--module, we obtain that any $\ueg$--subquotient occurring in
$[V({\bf P})^{\on {res}}_\ep]$ has highest weight less than or equal
to ${\bf P}$ with respect to our partial order. This proves the
surjectivity of $e$ and formula \eqref{multi}
\end{proof}

Given an irreducible $\uqg$--module $V$, we call the $\ueg$--module
$V^{\on{res}}_\ep$ the {\em specialization} of the module $V$.

\subsection{Finite-dimensional representations of $\fueg$}

A finite-dimensional representation $V$ of $\fueg$ is said to be of
type 1 if $k_i, i \in I$, act on $V$ semi-simply with eigenvalues in
$\ep_i^{\Z}$.

A vector $v$ in a $\fueg$--module $V$ is called a highest weight
vector if $x_{i,n}^+v=0$ for all $i \in I, n \in \Z$, and
\begin{equation}
\Phi_i^\pm(u) v = \Psi_i^\pm(u) v, \qquad \Psi_i^\pm(u)
\in \C[[u^{\pm 1}]] \qquad (i \in I).
\end{equation}
If in addition $V=\fueg\cdot v$, then $V$ is called a highest weight
representation with highest weight $(\Psi_i^\pm(u))_{i\in
I}$. Introduce the notation $\wt\Psi_i^+(u) = \Psi_i^+(u)/\Psi_i^+(0),
\wt\Psi_i^-(u) = \Psi_i^-(u)/\Psi_i^-(\infty)$ (these are the
eigenvalues of $k_i^{\mp 1} \Phi_i^\pm(u)$ on $v$). By \lemref{rel in
fin}, the elements $h_{i,ml_i}$ act on $V$ by $0$. Therefore
$\Psi_i^\pm(u)$ necessarily satisfies the property 
\begin{equation}    \label{proper}
\prod_{j=0}^{l_i-1} \wt\Psi_i^\pm(u\ep_i^{2j}) = 1 \qquad  (i\in I),
\end{equation}
cf. Section 6.5 of \cite{BK}. Note that the commuting elements
$h_{i,n}, n \neq ml_i$, are algebraically independent in $\ueg$, and
hence in $\fueg$. Therefore for any choice of $(\Psi_i^\pm(u))_{i\in
I}$ satisfying the condition \eqref{proper}, there exists a unique
irreducible highest weight representation with highest weight
$(\Psi_i^\pm(u))_{i \in I}$.

We will say that a polynomial $P(u) \in \C[u]$ is $l$--{\em acyclic}
if it is not divisible by $(1-au^l)$ (equivalently, the set of roots of
$P(u)$ does not contain a subset of the form $\{
a,a\ep^2,\ldots,a\ep^{2l-2} \}, a \in \C^\times$, where $\ep^2$ is a
primitive root of unity of order $l$).

The following statement essentially follows from the results of Beck
and Kac \cite{BK}.

\begin{thm}    \label{class fin}
Every irreducible finite-dimensional type 1 representation of $\fueg$ is a
highest weight representation. An irreducible highest weight
representation $V$ with highest weight $(\Psi_i^\pm(u))_{i \in I}$ is
finite-dimensional if and only if there exists an $I$--tuple of
polynomials with constant term $1$, ${\bf P}=(P_i(u))_{i\in I}$, where
$P_i(u)$ is $l_i$--acyclic, such that \bean \label{hw fin}
\Psi^\pm_i(u) = \gamma_i \ep_i^{\deg P_i}
\frac{P_i(u\ep_i^{-1})}{P_i(u\ep_i)}, \qquad i \in I, \eean where $\gamma_i = 1$, if
$\ep_i$ has odd order, $\gamma_i = \pm 1$, if $\ep_i$ has even order, and by
the rational function appearing in the right hand side we understand
its expansion in $u^{\pm 1}$.
\end{thm}

\begin{proof}
Theorem 6.3 of \cite{BK} describes the so-called diagonal
finite-dimensional representations of the non-restricted
specialization $\tueg$ (although $l$ is assumed in \cite{BK} to be
odd, the proof of Theorem 6.3 does not depend on this restriction). It
states in particular that these representations are highest weight
representations. Since we have a surjective homomorphism $\tueg
\twoheadrightarrow \fueg$, any irreducible $\fueg$--module gives rise
to an irreducible $\tueg$--module. 
By \lemref{rel in fin}, an irreducible $\tueg$--module obtained by
pull-back from an irreducible $\fueg$--module is diagonal.
Therefore we obtain that every
irreducible finite-dimensional $\fueg$--module is a highest weight
module.

Furthermore, according to Theorem 6.3 of \cite{BK}, the irreducible
highest weight representation of $\tueg$ with highest weight
$(\Psi_i^\pm(u))_{i \in I}$ is finite-dimensional if and only if
$\Psi^\pm_i(u)$ is a rational function $$f_i(u) =
\frac{P_i^{(1)}(u)}{P_i^{(2)}(u)}, \qquad i \in I,$$ where
$P_i^{(1)}(u), P_i^{(2)}(u)$ are two polynomials in $u$ of equal
degrees with non-zero constant coefficients; $f_i(u)$ is
regular at $0$ and at $\infty$, and $f_i(0) = f_i(\infty)^{-1}$.

In addition, the highest weight of an irreducible $\tueg$--module
obtained by pull-back from an irreducible $\fueg$--module of type 1
satisfies formula \eqref{proper} and the condition $\Psi^+_i(0) \in
\ep_i^{\Z}$. But then for each $i \in I$ there exists a unique acyclic
polynomial $P_i(u)$, such that formula \eqref{hw fin} holds. This
completes the proof.
\end{proof}

\begin{prop}    \label{fin by restr}
Let ${\bf P}=(P_i(u))_{i\in I}$ be an $I$--tuple of polynomials with
constant term $1$, such that $P_i(u)$ is $l_i$--acyclic for each $i
\in I$. Then the irreducible $\ueg$--module $V({\bf P})$ remains
irreducible when restricted to $\fueg$.
\end{prop}

\begin{proof}
This statement is proved in Theorem 9.2 of \cite{CP:root} in the case
when the order of $\ep$ is odd (although the algebra $\fueg$ in
\cite{CP:root} includes the elements $P_{i,n}$, the same proof works
for the algebra without the elements $P_{i,n}$, as in our definition
of $\fueg$). For general root of unity $\ep$, one repeats the proof of
\cite{CP:root} replacing $l$ by $l_i$, where appropriate.
\end{proof}

\begin{remark}    \label{get rid}
Let $I_{\on{ev}} = \{ i \in I| \ep_i \text{ has even order} \}$. The
algebra $\fueg$ has $2^{\# I_{\on{ev}}}$ non-isomorphic
one-dimensional type 1 representations, on which the $e_i$'s and
$f_i$'s act by zero, $k_i \in I\backslash I_{\on{ev}}$ act
identitically, and $k_i \in I_{\on{ev}}$ act by multiplication by $\pm
1$ (indeed, if $\ep_i$ is a root of unity of an even order $2l_i$,
then $\ep_i^{l_i} = -1$). By \thmref{class fin}, any irreducible
$\fueg$--module of type 1 is isomorphic to the tensor product of an
irreducible module with the highest weight of the form \eqref{hw fin},
where $\gamma_i=1$ for all $i \in I$, and a one-dimensional module.
Note that in contrast to $\fueg$, the algebra $\ueg$ has only one
one-dimensional type 1 representation for any $\ep$ (namely, the
trivial representation).

\propref{fin by restr} implies that any irreducible $\fueg$--module
may be extended to a $\ueg$--module, although this $\ueg$--module may
not be of type 1. Furthermore, an irreducible $\fueg$--module can be
extended to a $\ueg$--module in a {\em unique} way. This follows from
the Decomposition Theorem \ref{decomp thm} (which is proved in this
paper under the assumption that $(l,r^\vee)=1)$.
\end{remark}

\section{The $\ep$--characters}

\subsection{Definition}
Let $V$ be a finite-dimensional representation of $\ueg$. We have a
commutative subalgebra generated by the elements $k_i$, $P^\pm_{i,n}$  
($i \in I$, $n \in \Z$) acting on $V$. Let $v_1,\dots,v_d$
be a basis of common generalized eigenvectors of these elements. In
what follows we look at the corresponding generalized eigenvalues.

\begin{lem}\label{P act}
The common eigenvalues of $\mc P^+_i(u)$ are rational functions of the
form \be \Gamma_{i,n}^+(u) =
\frac{\prod_{j=1}^{j_{in}}(1-a_{inj}u)}{\prod_{m=1}^{m_{in}}(1-b_{inm}u)},
\ee where $j_{in},m_{in}$ are non-negative integers and
$a_{inj},b_{inm}$ are non-zero complex numbers. Moreover, the
generalized eigenvalues of $\mc P^-_i(u)$ are \be
\Gamma_{i,n}^-(u)=\frac{\prod_{j=1}^{j_{in}}
(1-(a_{inj}u)^{-1})}{\prod_{m=1}^{m_{in}}(1-(b_{inm}u)^{-1})}, \ee and
we have \be k_iv_n=q^{j_{in}-m_{in}}v_n, \qquad
\bin[k_i;l_i]v_n=\left[\begin{matrix} {j_{in}-m_{in}}\\{
l_i}\end{matrix}\right]_{\ep_i} v_n.  \ee
\end{lem}

\begin{proof}
In the case of $\uqg$ the analogous statement follows from Proposition
1 in \cite{FR}. Lemma \ref{P act} then follows from Proposition
\ref{spec rep}.
\end{proof}

The $q$--character of a $\uqg$--module $V$ defined in \cite{FR} (see
also \cite{FM}) encodes the generalized eigenvalues of $\mc P_i^+(u)$
on $V$. It is denoted by $\chi_q(V)$ and it takes values in the
polynomial ring $\Z[Y_{i,a}^{\pm 1}]_{i\in I}^{a\in \C^\times}$.

Similarly, we define the $\ep$--character of a finite-dimensional
$\ueg$--module $V$, denoted by $\chi_\ep(V)$, 
to be the element of the ring $\Z[Y_{i,a}^{\pm
1}]_{i\in I}^{a\in \C^\times}$ equal to \be
\chi_\ep(V)=\sum_{n=1}^d\prod_{i\in I}\left( \prod_{j=1}^{j_{in}}
Y_{i,a_{inj}}\prod_{m=1}^{m_{in}}Y_{i,b_{inm}}^{-1}\right).  \ee

\begin{thm}\label{spec char}
The map \bea \chi_\ep:\; \on{Rep}\ueg&\to &\Z[Y_{i,a}^{\pm 1}]_{i\in
I}^{a\in \C^\times},\\ V&\mapsto& \chi_\ep(V) \eea is an injective
homomorphism of rings. Moreover, for any irreducible
finite-dimensional $\uqg$--module $V$, $\chi_\ep(V^{\on{res}}_\ep)$ is
obtained from $\chi_q(V)$ by setting $q$ equal to $\ep$.
\end{thm}

\begin{proof}
The first part follows from Theorem 3 in \cite{FR} and \propref{spec
rep}. The second part follows from \propref{spec rep}.
\end{proof}

In particular, we obtain that the ring $\on{Rep}\ueg$ is commutative.

A monomial $m\in\Z[Y_{i,a}^{\pm 1}]_{i\in I}^{a\in \C^\times}$ is
called dominant if it does not contain factors $Y_{i,a}$ in negative
powers. The monomial in $\chi_\ep(V)$ corresponding to the highest
weight vector is always dominant. We call it the highest weight
monomial. The $i$-th fundamental representation $V_{\om_i}(a)$ is by
definition the irreducible $\ueg$--module whose highest weight monomial
is $Y_{i,a}$.

Theorems \ref{classe} and \ref{spec char} imply that the map \bea
\on{Rep}\ueg&\to& \Z[t_{i,a}]_{i\in I}^{a\in \C^\times}, \\
V_{\om_i}(a)&\mapsto & t_{i,a} \eea is an isomorphism of rings
(cf. the analogous statement in the case of $\on{Rep}\uqg$ in
\cite{FR}, Corollary 2).

\subsection{Properties of the $\ep$--characters}    \label{properties}
Using the definition, \propref{spec rep} and Theorem \ref{spec char},
we obtain that the $\ep$--characters satisfy combinatorial properties
similar to $q$--characters. In this section we list these
properties.

We have the ordinary character homomorphism $\chi: \on{Rep} \uegg \arr
\Z[y_i^{\pm 1}]_{i\in I}$: if $V = \oplus_\mu V_\mu$ is the weight
decomposition of $V$, then $\chi(V) = \sum_\mu \dim V_\mu \cdot
y^\mu$, where for $\mu = \sum_{i\in I} \mu_i \om_i$ we set $y^\mu =
\prod_{i\in I} y_i^{\mu_i}$.

Define the ring homomorphism $\beta:\;\Z[Y_{i,a}^{\pm 1}]_{i\in I}^{a\in
\C^\times}\to \Z[y_{i}^{\pm 1}]_{i\in I}^{a\in \C^\times}$
by the rule $\beta(Y_{i,a}^{\pm1})=y_i^{\pm 1}$. We will say that a
monomial $m \in \Z[Y_{i,a}^{\pm 1}]_{i\in I}^{a \in \C^\times}$ has
weight $\mu$ if $\beta(m) = y^\mu$. 

\begin{lem}\label{to finite}
The following diagram is commutative:
\be
\begin{CD} \on{Rep} \ueg @>{\chi_\ep}>> \Z[Y_{i,a}^{\pm 1}]_{i\in I}^{a
\in \C^\times}\\ @VV{\on{res}}V @VV{Y_{i,a}\mapsto y_i}V\\ \on{Rep}
\uegg @>{\chi}>> \Z[y_i^{\pm 1}]_{i\in I}
\end{CD}
\ee
\end{lem}

Given a subset $J$ of $I$, denote by $\ueg_J$ the subalgebra of $\ueg$
generated by $k_i$, $\bin[k_i;l_i]$, $P_{i,r}$,
$(x^\pm_{i,r})^{(n)}$, $i\in J$, $r\in\Z$, $n\in\Z_{\geq 0}$. Let
$\beta_J:\;\Z[Y_{i,a}^{\pm 1}]_{i\in I}^{a\in \C^\times}\to
\Z[Y_{i,a}^{\pm 1}]_{i\in J}^{a\in \C^\times}$ be a ring homomorphism
defined by the rule $\beta_J(Y_{i,a}^{\pm 1})=Y_{i,a}^{\pm 1}$ if
$i\in J$ and $\beta_J(Y_{i,a}^{\pm 1})=1$ if $i\not\in J$. Then we
have:

\begin{lem}    \label{comm diag for restr}
The following diagram is commutative.
\be
\begin{CD} \on{Rep} \ueg @>{\chi_\ep}>> \Z[Y_{i,a}^{\pm 1}]_{i\in I}^{a
\in \C^\times}\\ @VV{\on{res}}V @VV{Y_{i,a}\mapsto 1,\;\;i\not\in
J}V\\ \on{Rep} \ueg_J
@>{\chi_{\ep}^J}>> \Z[Y_{i,a}^{\pm 1}]_{i\in J}^{a
\in \C^\times}
\end{CD}
\ee
\end{lem}
Define  $A_{i,a} \in \Z[Y_{j,b}^{\pm 1}]_{j\in I}^{b\in\C^\times}$ by
the formula
\begin{equation}    \label{express}
A_{i,a}=Y_{i,a\ep_i^{-1}}\left(\prod_{j:\;j\neq i, \; C_{ji}\neq 0}
A_{i,j;a}\right)Y_{i,a\ep_i}, \quad
A_{i,j;a} = \left\{\begin{matrix}
Y_{j,a}^{-1} & {\rm if}\;\; C_{ji}=-1,\\
 Y_{j,a\ep^{-1}}^{-1}Y_{j,a\ep}^{-1} 
& {\rm if}\;\; C_{ji}=-2,\\
 Y_{j,a\ep^{-2}}^{-1}Y_{j,a}^{-1}  Y_{j,a\ep^{2}}^{-1}
 & {\rm if}\;\; C_{ji}=-3.
\end{matrix}\right.
\end{equation}
Thus $\beta(A_{i,a})$ corresponds to the simple root $\al_i$.

Using Theorem 4.1 in \cite{FM} and \propref{spec rep}, we obtain:

\begin{prop}
Let $V$ be an irreducible $\ueg$ module. The $\ep$--character of $V$ has
the form
\be
\chi_\ep(V)=m_+(1+\sum_p a_pM_p),
\ee
where $m_+$ is the highest weight monomial and for each $p$, $M_p$ is
a product of $A_{i,c}^{-1}$, $i\in I$, $c\in\C^\times$, $a_p$ is
a nonnegative integer.
\end{prop}

Define the rings $\mc K_i$, $i\in I$ by the formula
\be
{\mc K}_i = \Z[Y_{j,a}^{\pm 1}]_{j\in I;j\neq i}^{a \in
\C^\times} \otimes \Z[Y_{i,b} + Y_{i,b}
A_{i,b\ep_i}^{-1}]_{b \in \C^\times}.
\ee
Set
\be
\mc K=\bigcap\limits_{i\in I} \mc K_i\subset \Z[Y_{j,a}^{\pm 1}]_{j\in
I}^{a \in \C^\times}.
\ee
Corollary 5.7 in \cite{FM} and \propref{spec rep} imply:
\begin{prop}
The image of the $\ep$--character homomorphism equals
$\mc K$. Therefore
\be
\chi_\ep:\;\on{Rep}\ueg\to \mc K
\ee
is a ring isomorphism.
\end{prop}

\subsection{The $\ep$--characters of the small quantum group $\fueg$}

The $\ep$--characters of finite-dimen\-sional $\fueg$--modules are
defined similarly to the $\ep$--characters of $\ueg$--modules. Namely,
we consider the generalized eigenvalues $\Psi_i^\pm(u)$ of the series
$\Phi_i^{\pm}(u)$ on a given finite-dimensional $\fueg$--module
$V$. By Proposition \ref{fin by restr} and \remref{get rid}, any
irreducible $\fueg$--modules may be represented as the tensor product
of the restriction of a type one $\ueg$--module and a one-dimensional
sign representation.  Therefore we find from \lemref{P act} that the
generalized eigenvalues $\Psi_i^\pm(u)$ have the form \Ref{hw fin},
where $P_i(u) = \prod_{j=1}^{j_i} (1-a_{ij}u)$ is an $l_i$--acyclic
polynomial for each $i \in I$, and for each $i \in I$ the sign $\gamma_i$ in
\Ref{hw fin} is the same for all generalized eigenvectors. We attach to
a set of common generalized eigenvalues of this form, the monomial
$\prod_{i \in I} \prod_{j=1}^{j_i} Y_{i,a_{ij}}$.

Let $\Rep \fueg$ be the quotient of the Grothendieck ring of $\fueg$
by the ideal generated by the elements of the form $[V] - 1$, where
$[V]$ is any non-trivial one-dimensional $\fueg$--module. Introduce
the notation \be {\bf
Y}_{i,a}:=\prod_{j=0}^{l_i-1}Y_{i,a\ep_i^{2j}}. \ee Then we obtain an
injective ring homomorphism \be \chi^{\on{fin}}_\ep:\;\Rep \fueg\to
\C[Y_{i,a}]_{i\in I}^{a\in\C^\times}/({\bf Y}_{i,b}-1)_{i\in
I}^{b\in\C^\times}. \ee

Let $\pi$ be the ring homomorphism
\bea
\pi:\; \C[Y_{i,a}^{\pm1}]_{i\in I}^{a\in\C^\times}&\to&
\C[Y_{i,a}]_{i\in I}^{a\in\C^\times}/({\bf Y}_{i,b}-1)_{i\in
I}^{b\in\C^\times}, \\
Y_{i,a}&\mapsto& Y_{i,a}.
\eea
Then we have the following
\begin{lem}
Let $V$ be a finite dimensional $\ueg$ module with $\ep$--character
$\chi_\ep(V)$. Then the following diagram is commutative.
$\chi^{\on{fin}}_\ep\left( V|_{\fueg} \right) =
\pi(\chi_\ep(V))$.
\end{lem}

Since all irreducible $\fueg$--modules may be obtained by restriction
from irreducible $\ueg$--modules, we see that the study of
$\ep$--characters of $\fueg$--modules is reduced to the study of
$\ep$--characters of $\ueg$--modules.

\section{The Frobenius homomorphism}

In this section we introduce the Frobenius homomorphism following
Lusztig \cite{L}, \S 35. Because in the case of roots of unity of even
order this homomorphism is not discussed in detail in the literature
(except for \cite{L}), we first give a brief overview.

Recall that we denote by $\ep$ a fixed primitive root of unity of
order $s$. If $s$ is an odd number, relatively prime to $r^\vee$, then
we expect that there exists a Frobenius homomorphism $\ueg\to U\G$,
although we have not been able to locate such a homomorphism in the
literature (such a homomorphism has been constructed by Lusztig in
\cite{L1} in the case of $U_q\g$ with the above restrictions on $s$).
If $s$ is even, but not divisible by $4$ and by $r^\vee$ (so that $l$
is odd and not divisible by $r^\vee$), then presumably there exists a
homomorphism $\ueg\to U^{\on{res}}_{-1} \G$. If $s$ is even and the
above additional conditions are not satisfied, the construction
becomes more complicated.

The important aspect of Lusztig's definition of the Frobenius
homomorphism is that it uses the modified quantized enveloping algebra
$\dueg$ (it is denoted by $\dot{\mathbf U}$ in \cite{L}) instead of
$\ueg$. This appears to be necessary in order to give a uniform
definition of the Frobenius homomorphism for roots of unity of odd and
even orders (subject to some mild restrictions mentioned below). In
this section we apply Lusztig's definition of the Frobenius
homomorphism using the modified quantized enveloping algebra
$\dueg$. However, for the purposes of our paper the difference
between $\ueg$ and $\dueg$ is inessential, because the category of
finite-dimensional $\ueg$--modules of type 1 is equivalent to the
category of unital $\dueg$--modules.

\subsection{The modified quantized enveloping algebra}

Recall that we denote by $\Lambda$ the weight lattice of $\g$. For
$\la\in \Lambda$, \be U_{\la}=\uqg / \sum_{i\in
I}\uqg(k_i-q^{(\al_i,\la)}). \ee Here we set $(\al_0,\la) = -
(\al_{\on{max}},\la)$. Let $1_\la\in U_\la$ be the image of $1\in\uqg$
in $U_\la$. The space $U_\la$ is a left $\uqg$--module. For
$g\in\uqg$, we denote by $g1_\la$ the image of $g$ in $U_\la$.

Set \be \duqg=\oplus_{\la\in \Lambda}U_\la. \ee This is a
$\C(q)$--algebra with multiplication given by \bea 1_\la
1_\mu&=&\delta_{\la,\mu}1_\la,\\ x^\pm_i 1_\la&=&1_{\la\pm
\al_i} x^\pm_i.  \eea

Now we compare the algebra $\duqg$ with the modified quantized
enveloping algebra $\umod$ defined in \cite{L}. The latter is assigned
to a root datum, which consists of two lattices $X, Y$, a set $I$, a
pairing $\lan , \ran: Y \times X \to \Z$, a bilinear form $\cdot: \Z[I]
\times \Z[I] \to \Z$, and embeddings $I \hookrightarrow X, I
\hookrightarrow Y$ satisfying the conditions listed in \cite{L}, \S\S
1.1.1, 2.2.1.

We take as $I$ the set $\hat{I}$ of vertices of the Dynkin diagram of
$\G$, as $X$ the weight lattice $\Lambda$ of $\g$ (spanned by $\om_i, i \in
I$), and as $Y$ the coroot lattice $Q^\vee$ of $\g$ (spanned by
$\al_i^\vee, i \in I$). The pairing between $X$ and $Y$ is defined by
the formula $\lan \al_i^\vee,\om_j \ran = \delta_{ij}$ for all $i,j
\in I$. The bilinear form on $\Z[\hat{I}]$ is defined by the formula
$i \cdot j = (\al_i,\al_j)$ for $i,j \neq 0$, $0 \cdot j = -
(\al_{\on{max}},\al_j)$ for $j\neq 0$ and $0\cdot
0=(\al_{\on{max}},\al_{\on{max}})$.
The embedding $\hat{I} \hookrightarrow X$ is
defined by the formula $i \to \al_i, i \neq 0; 0 \to
-\al_{\on{max}}$. The embedding $\hat{I} \hookrightarrow Y$ is defined
by the formula $i \to \al_i^\vee, i \neq 0; 0 \to
-\al^\vee_{\on{max}}$. Attached to these data are associative algebras
over $\C(q)$, ${\mathbf U}$ and $\umod$, defined in \cite{L} (see \S
3.1.1, \S 23.1.1, and Corollary 33.1.5). Straightforward and explicit
comparison gives:

\begin{lem}
The algebra ${\mathbf U}$ is isomorphic to the algebra $\uqg$ extended
by the elements $k_i^{1/r_i}, i \in \hat{I}$. The algebra $\umod$ is
isomorphic to the algebra $\duqg$.
\end{lem}

A representation of $\duqg$ is called unital if $\sum_{\la\in \Lambda}
1_\la$ acts on it as the identity (note that the infinite sum
$\sum_\la1_\la$ is a well-defined operator on any $\duqg$--module). We
will consider only finite-dimensional unital $\duqg$--modules. These
are the finite-dimensional $\duqg$--modules, which do not contain
subspaces of positive dimension on which $\duqg$ acts by $0$.

Let $V$ be a finite-dimensional unital representation of $\duqg$. Then
using the projectors $1_\la$, we obtain a decomposition
$V=\oplus_{\la\in \Lambda} V_\la$. We define a $\uqg$--module structure on
$V$ by the rule $g\cdot v=g1_\la\cdot v$ for $g\in\uqg$, $v\in V_\la$.

Conversely, if $V$ is a type 1 finite-dimensional $\uqg$--module, then
we have the weight decomposition $V=\oplus_\la V_\la$. Define the
action of $\duqg$ in $V$ by letting $1_\la$ act as the identity on
$V_\la$ and by zero on $V_\mu$, $\mu\neq \la$.

Hence we obtain the following result (cf. \cite{L}, \S 23.1.4):

\begin{lem}\label{dot cat} The category of finite-dimensional
$\uqg$--modules of type 1 is equivalent to the category of
finite-dimensional unital $\duqg$--modules.
\end{lem}

Note that strictly speaking $\duqg$ is not a Hopf algebra. However, we
define tensor product of $\duqg$--modules using Lemma \ref{dot cat},
see also \cite{L}, \S 23.1.5.

The automorphisms $T_i$ on $\uqg$ induce automorphisms of
$\duqg$, which we also denote by $T_i$.

\subsection{Specialization to the root of unity}

Next we define $\duqgr$ as the $\C[q,q^{-1}]$--subalgebra of $\duqg$,
generated by the elements $(x^\pm_i)^{(n)}1_\la, i \in \hat{I}, n \geq
0, \la \in \Lambda$. Note that we have the following relations in $\duqgr$,
$$
(x^\pm_i)^{(n)}1_\la=1_{\la\pm n\al_i}(x^\pm_i)^{(n)}.
$$
Then we set $\dueg = \duqgr \otimes_{\C[q,q^{-1}]} \C_\ep$.

\begin{lem}\label{incl}
The algebras $\duqgr$ and $\dueg$ are generated by
$(x_{i,r}^{\pm})^{(n)}1_\la, i \in I, r \in \Z, n \geq 0, \la \in \Lambda$.

If $J\subset I$, then the maps 
$$\duqgr_J \to \duqgr, \qquad \dueg_J\to\dueg,$$ 
$$(x_{i,r}^{\pm})^{(n)}1_\la \mapsto
(x_{i,r}^{\pm})^{(n)}1_\la, \qquad i\in J,\; r\in \Z,\; n\in\Z_{>0}$$
extend to injective algebra homomorphisms. Here we extended a weight
$\la$ of $\G_J$ to a weight of $\G$ by the rule $\lan \al_i^\vee,\la
\ran=0$ if $i\in I\setminus J$.
\end{lem}

\begin{proof}
This lemma follows from \thmref{trdec} and the commutation relations
in $\uqg$.
\end{proof}

\subsection{Definition of the Frobenius homomorphism}
Recall that $\ep^s=1$, $l=s$ if $s$ is odd and $l=s/2$ if $s$ is
even. 

From now on, we impose the following restrictions on $s$. First we
assume that $s>2$ (if $s=1,2$ then the Frobenius homomorphism is just
the identity map). Next, in this section we exclude from consideration
the case $\G=\widehat{\sw}_{2n+1}$ because of condition 35.1.2 b) in
\cite{L}. We will deal with it in Section \ref{odd sl}.  Finally, due
to condition 35.1.2 a) in \cite{L} we assume that $l\neq 2,4$ in the
case $r^\vee=2$, and $l\neq 2,3,6$ in the case $r^\vee=3$.

In what follows a weight $\la \in \Lambda$ is called {\em $l$--admissible}
if $\lan \al_i^\vee,\la \ran$ is divisible by $l_i$ for all $l$. In
other words, $\la$ is $l$--admissible if $\la=\sum_{i\in I} a_i\om_i$,
where each $a_i$ is divisible by $l_i$. Denote the lattice of all
$l$--admissible weights by $\Lambda_l$, and for each $\la = \sum_{i\in I}
a_i\om_i\in \Lambda_l$, denote by $\la/{\bf l}$ the weight $\sum_{i\in I}
(a_i/l_i) \om_i$.

In \cite{L}, \S 35.1.6 an algebra $\umod{}^*$ over $\C[q,q^{-1}]$ is
introduced. It is defined in the same way as $\umod{}$ with respect to
the Cartan datum dual to the Cartan datum used in the definition of
$\umod$ (see \cite{L}, \S 2.2.5). Let $R$ be a $\C[q,q^{-1}]$--module,
on which the operator of multiplication by $q$ has order $s$. Set
$_R\umod = \umod \otimes_{\C[q,q^{-1}]} R, {}_R\umod = \umod{}^*
\otimes_{\C[q,q^{-1}]} R$. Then by Theorem 35.1.9 in \cite{L}, there
exists an algebra homomorphism ${}_R\umod \to {}_R\umod{}^*$.

In the case at hand this result translates as follows. If $l$ is
relatively prime to $r^\vee$, define $\duqges$ to be $\duegs$, where
$\ep^* = \ep^{l^2}$. If $l$ is divisible by $r^\vee$, define $\duqges$
to be $\dot{U}^{\on{res}}_{\ep^*} \G^L$, where $\ep^* =
\ep^{l^2/r^\vee}$. Here $\G^L$ is the twisted affine Kac-Moody algebra,
whose Cartan matrix is the transpose of the Cartan matrix of $\G$. The
algebra $\dot{U}^{\on{res}}_{\ep^*} \G^L$ is defined in the same way
as $\dueg$ above (using the Drinfeld-Jimbo realization, with the generators
labeled in a way compatible with the labeling of the generators of
$\dueg$).

Then there exists an algebra homomorphism \bea
\on{Fr}:\;\dueg&\to&\duqges,\\ (x_i^{\pm})^{(m)}1_\la&\mapsto& (\bar
x_i^{\pm})^{(m/l_i)}\bar 1_{\la/{\bf l}} , \qquad m/l_i \in\Z, \la\in
\Lambda_l,\\ (x_i^{\pm})^{(m)}1_\la&\mapsto& 0,\qquad \qquad\qquad \;
\; \quad {\rm otherwise}. \eea Here and below we put a bar over the
elements of the target algebra to avoid confusion.

The homomorphism $\on{Fr}$ is called the {\em quantum Frobenius
homomorphism}.

\medskip

\begin{remark}
In \cite{L}, the algebra $\uqg$ is defined over ${\mathbb Q}(q)$ and
the algebra $\uqgr$ is defined over ${\mathbb Z}[q,q^{-1}]$. However,
all of the results of \cite{L} referred to above, remain valid if we
replace ${\mathbb Q}(q)$ by $\C(q)$ and ${\mathbb Z}[q,q^{-1}]$ by
$\C[q,q^{-1}]$.

The defining relations of the algebra $\umod$ used in \cite{L} are
strictly speaking different from those used here. Namely, the
relations between the generators $x^+_i, i \in \hat{I}$ (and $x^-_i, i
\in \hat{I}$) are described in \cite{L} in terms of a certain bilinear
form (see \cite{L}, \S\S 1.2.3, 1.2.4) instead of the quantum Serre
relations that we use here. However, it is known that the quantum
Serre relations are included into these relations (see \cite{L}, \S
1.4), and therefore there is a surjective homomorphism $\duqgr \to
\umod$. When $q=1$, this homomorphism is actually an isomorphism by a
theorem of Gabber--Kac, and the same is true for $q=-1$ by \cite{L},
\S 33.2 (unless $\G=\widehat{\sw}_{2n+1}$ which we have excluded from
consideration in this subsection). The homomorphism $\on{Fr}$ is the
composition of the homomorphism $\duqgr \to \umod$ with the Frobenius
homomorphism from \cite{L}, Theorem 35.1.9 (when $R=\C_\ep$).\qed

\end{remark}

\medskip

In order to avoid twisted quantum affine algebras, we assume from now
on that {\em $l$ is relatively prime to $r^\vee$}. We will describe
the $q$--characters of twisted quantum affine algebras and the
Frobenius homomorphism involving them in a separate paper.

Since $l$ is relatively prime to $r^\vee$, we have: $\ep_i=\ep$ and
$l_i=l$ for all $i \in \hat{I}$. Therefore $\Lambda_l = l \cdot
\Lambda$ and $\la/{\bf l} = \la/l$ for $\la \in \Lambda_l$.

Note that $\ep^*=\ep^{l^2}$ is equal to $1$ if $l$ is odd and $s=l$,
or if $l$ is even, and to $-1$ if $l$ is odd and $s=2l$. Note also that
$-\ep^*=(-\ep)^l$. 
By Proposition 33.2.3 in \cite{L}, the algebras ${\dot{U}^{\on{res}}_1
\G}$ and ${\dot{U}^{\on{res}}_{-1} \G}$ are isomorphic if
$\G\neq \widehat{\sw}_{2n+1}$.

In the same way as in the proof of \lemref{dot cat} we obtain the
following result (see the beginning of \secref{type 1} for the
definition of $\ueg$--modules of type 1).

\begin{lem}
The category of finite-dimensional $\ueg$--modules of type 1 is
equivalent to the category of finite-dimensional unital
$\dueg$--modules.
\end{lem}

Therefore, the quantum Frobenius homomorphism induces a map of
Grothendieck rings \be \on{Fr}^*:\;\on{Rep} \uegs \to \on{Rep}
\ueg. \ee By \cite{L}, \S 35.1.10, this map is a ring homomorphism.

\begin{lem}\label{frob curr}
We have $\on{Fr}((x^{\pm}_{i,r})^{(m)}1_\la)=0$ if $m$ is not
divisible by $l$ or if $\la\not\in \Lambda_l$. For $\la\in \Lambda_l$, \be
\on{Fr}((x^{\pm}_{i,r})^{(nl)}1_\la)=\left\{
\begin{matrix} 
(\bar x^{\pm}_{i,r})^{(n)}\bar 1_{\la/l}, & l {\rm \;\; is \;\; odd}\\
o(i)^{nr}(\bar x^{\pm}_{i,r})^{(n)}\bar 1_{\la/l}, & l {\rm \;\; is
\;\; even.}
\end{matrix}\right.
\ee
\end{lem}
\begin{proof}
We use the fact that the automorphisms $T_{i}, i \in \hat{I}$, commute
with the Frobenius homomorphism (see \cite{L}, \S 41.1.9), and so do
the automorphisms $\sigma$ of the Dynkin diagram of $\G$ (see
\secref{qaal}). Hence the automorphisms $T_{\check\om_i}, i \in I$,
also commute with the Frobenius homomorphism. From this we find: \bea
\on{Fr}((x^{\pm}_{i,r})^{(m)}1_\la)=\on{Fr}(o(i)^{mr}
T_{\check\om_i}^{\mp r}(x^{\pm}_{i})^{(m)}1_{\check\om_i^{\pm
r}(\la)})=\\ = o(i)^{mr} \bar{T}_{\check\om_i}^{\mp
r}\on{Fr}((x^{\pm}_{i})^{(m)}1_{\check\om_i^{\pm r}(\la)}).  \eea So,
it is zero if $m$ is not divisible by $l$ or if $\la\not\in
\Lambda_l$.  If $m=nl$ we obtain \bea
\on{Fr}((x^{\pm}_{i,r})^{(m)}1_\la)=
o(i)^{mr}\bar{T}_{\check\om_i}^{\mp r}(\bar
x^\pm_{i})^{(n)}\bar 1_{\check\om_i^{\pm r}(\la)/l}=\\ =
o(i)^{mr+nr}(\bar x^{\pm}_{i,r})^{(n)}\bar 1_{\la/l}.  \eea
If $l$ is odd then $m$ and $n$ have the same parity, and all the signs
cancel. If $l$ is even, then $m$ is also even, and we acquire the sign
$o(i)^{nr}$.
\end{proof}

\subsection{The case $\G=\widehat{\sw}_{2n+1}$}\label{odd sl}
Let now $\G=\widehat{\sw}_{2n+1}$.  We define the Frobenius
homomorphism $\on{Fr}: \dot{U}_{\ep^*} \widehat{\sw}_{2n+1}
\to \dot{U}_{\ep} \widehat{\sw}_{2n+1}$ using the
generators $(x^\pm_{i,r})^{(m)}$ as in Lemma \ref{frob curr}: \bea
\on{Fr}((x^{\pm}_{i,r})^{(nl)}1_\la)&=&\left\{
\begin{matrix} 
(\bar x^{\pm}_{i,r})^{(n)}\bar 1_{\la/l} & l {\rm \;\; is \;\; odd}, \;\;\la\in
\Lambda_l,\\
o(i)^{nr}(\bar x^{\pm}_{i,r})^{(n)}\bar 1_{\la/l} & l {\rm \;\; is
\;\; even},\;\;\la\in \Lambda_l
\end{matrix}\right.
\eea
and $\on{Fr}((x^{\pm}_{i,r})^{(m)}1_\la)=0$ if $m$ is not divisible by $l$
or $\la\not\in \Lambda_l$.

\begin{lem}
These formulas give rise to a well-defined homomorphism of algebras.
Moreover, the induced map $\on{Fr}^*:\; \on{Rep}{U}^{\on{res}}_{\ep^*}\G\to
\on{Rep}{U}^{\on{res}}_{\ep}\G$ is a ring homomorphism.
\end{lem}

\begin{proof}
We embed the Dynkin diagram of $\widehat{\sw}_{2n+1}$ into the Dynkin
diagram of $\widehat{\sw}_{2n+2}$ in such a way that the numbers
$o(i)$ coincide. By Lemma \ref{incl}, we have the corresponding
embeddings $\dot{U}^{\on{res}}_\ep\widehat{\sw}_{2n+1}\to
{\overset{\cdot} U}_\ep\widehat{\sw}_{2n+2}$ and
$\dot{U}_{\ep^*}\widehat{\sw}_{2n+1}\to
{\dot{U}}_{\ep^*}\widehat{\sw}_{2n+2}$. It follows from
\lemref{frob curr} that the image of
$\dot{U}^{\on{res}}_\ep\widehat{\sw}_{2n+1}$ under the
Frobenius homomorphism ${\dot{U}}_\ep\widehat{\sw}_{2n+2}
\to {\dot{U}}_{\ep^*}\widehat{\sw}_{2n+2}$ is contained in
${\dot{U}}_{\ep^*}\widehat{\sw}_{2n+1}$. Hence we obtain an
algebra homomorphism
$\dot{U}^{\on{res}}_\ep\widehat{\sw}_{2n+1} \to
{\dot{U}}_{\ep^*}\widehat{\sw}_{2n+1}$, given by the above
formula. The lemma follows.
\end{proof}

By a direct computation similar to the one used in the proof of Lemma
\ref{frob curr} we obtain the following formulas for the Frobenius
homomorphism in the case of $U^{\on{res}}_\ep\widehat{\sw}_{2n+1}$ in
terms of the Drinfeld-Jimbo generators.

\begin{lem} Let $\G=\widehat{\sw}_{2n-1}$, then
\bea \on{Fr}((x_i^{\pm})^{(m)}1_\la)&=& (\bar
x_i^{\pm})^{(m/l)}\bar 1_{\la/l} , \qquad \qquad\;\;\;\; m/l\in\Z, \la\in
\Lambda_l, i\in \hat I,\\ 
\on{Fr}((x_i^{\pm})^{(m)}\bar 1_\la)&=& 0,\qquad\;\;\;\;\qquad
\qquad\qquad \; \; \; {\rm otherwise}.  \eea
\end{lem}

\section{Frobenius homomorphism and $\ep$--characters}

\subsection{Properties of the Frobenius homomorphism}

Recall the automorphism $\tau_a$ defined by formula \eqref{tau}. It
gives rise to an automorphism of $\dueg$ in an obvious way. Using
Lemma \ref{frob curr}, we obtain:

\begin{lem}\label{shift inv}
We have: $\on{Fr}\circ\tau_a=\tau_{a^l}\circ
\on{Fr}$. In particular, $\on{Fr}\circ\tau_{\ep^2}=\on{Fr}$.
\end{lem}

Denote by $\tau^*_a$ the endomorphism of $\on{Rep} \ueg$ induced by
$\tau_a$. The next lemma follows from Lemmas \ref{incl} and \ref{frob
curr}.

\begin{lem}\label{restr}
Let $J\subset I$ correspond to a Dynkin subgraph of $\g$ and let $g_J$
be the corresponding Lie algebra. If the functions
$o(i)$ for diagrams $I$ and $J$ coincide or $l$ is odd, then 
the following diagram is commutative.
\be
\begin{CD} \on{Rep}{U}^{\on{res}}_{\ep^*}\G @>{\on{Fr}^*}>>\on{Rep}\ueg\\
@VV{\on{res}}V @VV{\on{res}}V\\ \on{Rep}{U}^{\on{res}}_{\ep^*}\G_J
@>{\on{Fr}^*_J}>> \on{Rep}\ueg_J
\end{CD}
\ee
If the functions
$o(i)$ for diagrams $I$ and $J$ do not coincide and $l$ is even then
${\on {Fr}}^*\circ \bar \tau_{-1}^*=\tau_\ep\circ {\on {Fr}}^*$ and
the following diagram is commutative.
\be
\begin{CD} \on{Rep}{U}^{\on{res}}_{\ep^*}\G @>{\on{Fr}^*}>>\on{Rep}\ueg\\
@VV{\on{res}}V @VV{\on{res}}V\\ \on{Rep}{U}^{\on{res}}_{\ep^*}\G_J
@>{\on{Fr}^* \circ \bar \tau^*_{-1}}>> \on{Rep}\ueg_J
\end{CD}
\ee
\end{lem}

The algebras $\dot{U}^{\on{res}}_\ep\g$,
$\dot{U}^{\on{res}}_{\ep^*}\g$ 
and the corresponding Frobenius map are defined similarly to the
affine case. Then the following property 
follows immediately from the definitions in terms of
the Drinfeld-Jimbo generators.

\begin{lem}    \label{rest to finite}
The following diagram is commutative:
\be
\begin{CD} \on{Rep}{U}^{\on{res}}_{\ep^*}\G @>{\on{Fr}^*}>>\on{Rep}\ueg\\
@VV{\on{res}}V @VV{\on{res}}V\\ \on{Rep}{U}^{\on{res}}_{\ep^*}\g
@>{\on{Fr}^*}>> \on{Rep}\uegg
\end{CD}
\ee
\end{lem}

\subsection{The decomposition theorem}

For a polynomial $P(u)$ with constant term 1, define polynomials with
constant term 1, $P^0(u)$ and $P^1(u)$, to be such that
$P(u)=P^0(u)P^1(u)$, $P^0(u)$ is not divisible by $1-au^l$, for any
$a\in\C^\times$ and $P^1(u)=R(u^l)$ for some polynomial $R$. If ${\bf
P} = (P_i(u))_{i\in I}$, we write: ${\bf P}^0 = (P^0_i(u))_{i\in I}$
and ${\bf P}^1 = (P^1_i(u))_{i\in I}$.

Let $V$ be a ${U}^{\on{res}}_{\ep^*}\G$ module. We call the
$\ueg$--module $\on{Fr}^*(V)$ obtained by pull-back of $V$ via the
quantum Frobenius homomorphism the {\em Frobenius pull-back} of $V$.

\begin{thm}\label{decomp thm}
Let $V({\bf P})$ be an irreducible $\ueg$ module with Drinfeld
polynomials ${\bf P}$. Then $V({\bf P})\simeq V({\bf P}^0)\otimes
V({\bf P}^1)$. Moreover $V({\bf P}^1)$ is the Frobenius pull-back of
an irreducible $U^{\on{res}}_{\ep^*} \ghat$--module and $V({\bf P}^0)$
is irreducible over $U^{\on{fin}}_\ep \ghat$.
\end{thm}

\begin{proof}
This theorem is proved in \cite{CP:root}, Theorems 9.1--9.3, in the
case when $\ep$ is a root of unity of odd order (i.e., $s=l$ and $l$
is odd). The proof given in \cite{CP:root} goes through for other $l$
(under our assumption that $l$ is relatively prime to $r^\vee$) with
the following changes: Lemma 9.5 of \cite{CP:root} should be replaced
by Lemma \ref{frob curr} and the diagram (49) of \cite{CP:root} in the
case of even $l$ should be replaced by the identity $\on{Fr} \circ \,
ev_b=ev_{(\ep b)^l} \circ \on{Fr}$ (which is proved in the same way as
the commutativity of the diagram (49)). Finally, the proof uses the
description of the Drinfeld polynomials of the Frobenius pull-backs of
irreducible $U^{\on{res}}_{\ep^*}$--modules. These Drinfeld
polynomials are determined in \thmref{hw of frob} below.
\end{proof}

\thmref{decomp thm} and \propref{fin by restr} imply:

\begin{cor}\label{break cor}
There is an isomorphism of vector spaces
\be
\Rep \ueg \simeq \Rep U_{\ep^*}^{\on{res}} \G \otimes \Rep \fueg.
\ee
\end{cor}

The decomposition in Corollary \ref{break cor} is not a decomposition
of rings. However, the Frobenius map $\on{Fr}^*:\;\Rep
U_{\ep^*}^{\on{res}}\G\to \Rep \ueg$ is a ring homomorphism. Therefore
we have a natural ring structure on the quotient of $\Rep \ueg$ by its
ideal generated by the augmentation ideal of $\Rep
U_{\ep^*}^{\on{res}} \G$. By \corref{break cor}, this quotient is
isomorphic to $\Rep \fueg$ as a vector space. 
We call the induced multiplication on
$\Rep \fueg$ the {\em factorized tensor product}. It would be
interesting to extend it to the level of the category of
finite-dimensional representations of $\fueg$.

The statement of \thmref{decomp thm} remains true in the case when $l$
is relatively prime to $r^\vee$. However, in the course of proving it
we need to use information about the Frobenius homomorphism which in
that case takes values in the (modified) twisted quantum affine
algebra. This case will be discussed in a separate paper.

\subsection{The $\ep$--characters of the Frobenius pull-backs}

Recall the notation \be {\bf
Y}_{i,a}:=\prod_{j=0}^{l-1}Y_{i,a\ep^{2j}}. \ee Note that because $l$
is relatively prime to $r^\vee$, and hence by $r_i$, we have: ${\bf
Y}_{i,a}=\prod_{j=0}^{l-1}Y_{1,a\ep_i^{2j}}$. The monomial ${\bf
Y}_{i,a}$ corresponds to the polynomial $(1-a^lu^l)$.

Note that we have a homomorphism $\chi_{\ep^*}: \Rep
U^{\on{res}}_{\ep^*} \ghat \to \Z[\bar Y_{i,a}^{\pm1}]_{i\in
I}^{a\in\C^\times}$.

\begin{lem}\label{h w}
Let $V$ be an irreducible representation of
${U}^{\on{res}}_{\ep^*}\wh{\sw}_2$ with Drinfeld polynomial $\bar
P(u)$. Then the Drinfeld polynomial of $\on{Fr}^*(V)$ is $P(u) = \bar
P(u^l)$, if $l$ is odd, and $\bar P(o(1)u^l)$, if $l$ is even.

Moreover, $\chi_\ep(\on{Fr}^*(V))$ is obtained from $\chi_{\ep^*}(V)$
by replacing $\bar Y_{1,a^l}^{\pm 1}$ with ${\bf Y}_{1,a}^{\pm 1}$, if
$l$ is odd, and with ${\bf Y}_{1,a\ep^{\theta(1)}}^{\pm 1}$, where
$\theta(1) = (1-o(1))/2$, if $l$ is even.
\end{lem}

\begin{proof}
Note that if $a^l=b^l$ then $(a/b)^l=1$ and ${\bf Y}_{i,a}={\bf
Y}_{i,b}$. Therefore the rule described in the theorem is unambiguously
defined.

Recall that $$\on{Fr}^*: \on{Rep} {U}^{\on{res}}_{\ep^*}\wh{\sw}_2 \to
\on{Rep} {U}^{\on{res}}_{\ep}\wh{\sw}_2,$$ as well as $\chi_\ep$ and
$\chi_{\ep^*}$, are ring homomorphisms. Since the ring
${U}^{\on{res}}_{\ep^*}\wh{\sw}_2$ is generated by the fundamental
representation $V_{\om_1}(a)$, it is sufficient to prove the lemma
when $V=V_{\om_1}(a)$.

In this case the Drinfeld polynomial is $\bar P(u)=1-au$. Let $v$ be
the highest weight vector in $V$ and $\on{Fr}^*(V_{\omega_1}(a)))$. By
Lemma 5.1 in \cite{CP:root}, we obtain: \be (-1)^m\ep^{-m^2+m}k^m
P_{m}v=(x_{1,0}^+)^{(m)}(x_{1,1}^-)^{(m)}v. \ee
(Note that $P_m$ in this paper differs from $P_m$ in \cite{CP} 
by the factor $q^{-m}$.)
So, $P_mv=0$ unless $m$
is divisible by $l$.

If $m=nl$, then for odd $l$ we obtain: \be
(x_{1,0}^+)^{(m)}(x_{1,1}^-)^{(m)}v=(\bar x_{1,0}^+)^{(n)}(\bar
x_{1,1}^-)^{(n)}v = (-1)^n (\ep^{l^2})^{-n^2+n}\bar k^n\bar P_{n}v. \ee
We have: $\bar P_nv=0$, unless $n=1$, and $\bar P_1v=-av$. In
addition, $\bar k v = \ep^{l^2} v, k v = \ep^l v$. Hence we find: $P_m
= 0$, unless $m=l$, and $P_l v = -a v$.

For even $l$ we obtain: \be (x_{1,0}^+)^{(m)}(x_{1,1}^-)^{(m)} v =
o(1)^n (\bar x_{1,0}^+)^{(n)}(\bar x_{1,1}^-)^{(n)}v = (-1)^n o(1)^n
(\ep^{l^2})^{-n^2+n}\bar k^n\bar P_{n}v. \ee We have: $\bar P_n v=0$
unless $n=1$, and $\bar P_1v=-av$. Hence we find that $P_m = 0$,
unless $m=l$, and $P_l v =-o(1)a v$.

It remains to determine the $\ep$--character of
$\on{Fr}^*(V_{\omega_1}(a^l))$. Consider the case when $l$ is
odd. Since the Drinfeld polynomial of this module is $1-a^lu^l$, we
find that the highest weight monomial is ${\bf Y}_{1,a}$. By
\propref{spec rep} and \thmref{spec char}, the other monomial in it
must be obtained by specialization of a monomial in the $q$--character
of the irreducible $U_q \widehat{\sw}_2$--module with the highest
weight monomial $\prod_{j=0}^{l-1}Y_{1,aq^{2j}}$. Moreover, by
\lemref{to finite}, the degree of this monomial must be equal to
$-l$. There is only one monomial of such degree, i.e., the lowest
weight monomial, and it is equal to
$\prod_{j=0}^{l-1}Y^{-1}_{1,aq^{2j+2}}$. Therefore we obtain that the
$\ep$--character of $\on{Fr}^*(V_{\omega_1}(a^l))$ is equal to ${\bf
Y}_{i,a} + {\bf Y}_{i,a}^{-1}$. On the other hand, we have:
$\chi_{\ep^*}(V_{\omega_1}(a^l)) = \bar Y_{1,a^l} + \bar
Y_{1,a^l}^{-1}$, cf. Section \ref{classic},
and the second assertion of the lemma follows for odd
$l$.

If $l$ is even, we obtain by the same argument that the
$\ep$--character of $\on{Fr}^*(V_{\omega_1}(a^l))$ is equal to ${\bf
Y}_{i,a\ep^{\theta(1)}} + {\bf Y}_{i,a\ep^{\theta(1)}}^{-1}$. This
proves the second assertion of the lemma for even $l$.
\end{proof}

We can now obtain the description of the $\ep$--characters of
Frobenius pull-backs for general $\g$. Recall that we are under the
assumption that $l$ is relatively prime to $r^\vee$.

Set $\theta(i) = (1-o(i))/2$.

\begin{thm}    \label{hw of frob}
Let $V$ be an irreducible representation of
${U}^{\on{res}}_{\ep^*}\ghat$ with Drinfeld polynomials $P_i(u), i \in
I$. Then the $i$-th Drinfeld polynomial of $\on{Fr}^*(V)$ is equal to
$P_i(u^l)$, if $l$ is odd, and to $P_i(o(i)u^l)$ if $l$ is
even.

Moreover, if $l$ is odd, then $\chi_\ep(\on{Fr}^*(V))$ is obtained from
$\chi_{\ep^*}(V)$ by replacing $\bar Y_{i,a^l}^{\pm 1}$ with ${\bf
Y}_{i,a}^{\pm 1}$. If $l$ is even, then $\chi_\ep(\on{Fr}^*(V))$ is
obtained from $\chi_{\ep^*}(V)$ by replacing $\bar Y_{i,a^l}^{\pm 1}$
with ${\bf Y}_{i,a\ep^{\theta(i)}}^{\pm 1}$.
\end{thm}

\begin{proof}
Restricting to the subalgebras $U^{\on{res}}_{\ep_i}\widehat\sw_2
\subset \ueg$ and using Lemmas \ref{restr} and \ref{comm diag for
restr}, we obtain the statement of the theorem from \lemref{h w}.
\end{proof}

\thmref{hw of frob} allows us to write down explicit formulas for the
$\ep$--characters of Frobenius pull-backs from the $\ep^*$--characters
of irreducible ${U}^{\on{res}}_{\ep^*}\ghat$--modules, which we now
set out to determine.

\subsection{The $\ep^*$--characters of irreducible
${U}^{\on{res}}_{\ep^*}\ghat$--modules}\label{classic}

First, we consider the case when $\ep^*=1$ (i.e., $l$ is odd and
$s=l$, or $l$ is even). There is a surjective homomorphism
$U_1^{\on{res}} \G \to U \G$ (where $U \G$ is the universal enveloping
algebra of $\G=\g[t,t^{-1}]$), which sends the generators
$(x^\pm_{i})^{(n)}$ of $U_1^{\on{res}} \G$ to the corresponding
generators of $U \G$. In fact, this map gives rise to an isomorphism
between $U \G$ and the quotient of $U_1^{\on{res}} \G$ by the ideal
generated by the central elements $(K_i-1), i \in I$ (see \cite{CP},
Proposition 9.3.10). Hence the category of finite-dimensional type 1
$U_1^{\on{res}} \G$--modules is equivalent to the category of
finite-dimensional $\G$--modules, on which the Cartan subalgebra $\h
\subset \g \subset \G$ acts diagonally (we call them weight
modules).

The description of irreducible finite-dimensional representations of
$\G$ is as follows \cite{C,CP1}. Consider the ``evaluation
homomorphism'' $\phi_a: \G = \g[t,t^{-1}] \arr \g$ corresponding to
evaluating a Laurent polynomial in $t$ at a point $a \in
\C^\times$. For an irreducible $\g$--module $V_\la$ with highest
weight $\la$, let $V_\la(a)$ be its pull-back under $\phi_a$ to an
irreducible representation of $\G$. Then $V_{\la_1}(a_1) \otimes
\ldots \otimes V_{\la_n}(a_n)$ is irreducible if $a_i \neq a_j,
\forall i\neq j$, and these are all irreducible finite-dimensional
representations of $\G$ up to an isomorphism.

Therefore in order to obtain $\chi_1(V)$ for an arbitrary irreducible
${U}^{\on{res}}_{\ep^*}\ghat$--module, it suffices to know
$\chi_1(V_\la(a))$. Those is found by explicitly computing the
image of ${\mc P}^\pm_i(u)$ under the evaluation homomorphism. The
answer is the following. Let $\chi(V_\la)$ be the ordinary character of
the $\g$--module $V_\la$, considered as a polynomial in $y_i^{\pm 1},
i \in I$, as in \secref{properties}.

\begin{lem}
$\chi_1(V_\la(a))$ is obtained from $\chi(V_\la)$ by replacing each
$y_i^{\pm 1}$ with $\bar{Y}_{i,a}^{\pm 1}$.
\end{lem}

Combining this statement with \thmref{hw of frob}, we obtain a
complete description of the $\ep$--characters of irreducible Frobenius
pull-backs in the case when $\ep^*=1$.

Next, consider the case when $\ep^*=-1$ (i.e., $s=2l$ and $l$ is
odd). Introduce a function $\psi: I \to \{ \pm 1 \}$ by the rule:
$\psi(i) \neq \psi(j)$, whenever $i\neq j$,
$(\al_i,\al_j) \neq 0$, and $\on{min}
\{ r_i,r_j \} = 1$. We fix this function by the requirement that
$o(i)=\psi(i)$ if $r_i=1$.  


We say that a monomial $m \in \Z[\bar{Y}_{i,a}^{\pm 1}]_{i\in
I,a\in\C^\times}$ is supported at $a_0\in\C^\times$ if
$m\in\Z[\bar{Y}_{i,a_0\psi(i)}^{\pm1}]_{i\in I,k\in\Z}$. Any monomial
can be written in a unique way as a product of monomials supported at
some $a_1,\ldots,a_k \in \C^\times$, $a_i\neq a_j$.

\begin{lem}\label{minus 1 lem}
Let $V$ be an irreducible finite-dimensional $U_{-1}^{\on{res}}
\G$--module (of type 1), such that the highest weight monomial in
$\chi_{-1}(V)$ is the product $m_1 \ldots m_k$, where $m_i$ is
a monomial 
supported at $a_i \in \C^\times$ and has weight $\la_i$. Assume that
all $a_i$ are distinct, $a_i\neq a_j$. Then $V
\simeq V_1 \otimes \ldots \otimes V_k$, where $V_i$ is an irreducible
finite-dimensional $U_{-1}^{\on{res}} \G$--module with highest weight
monomial $m_i$.

Moreover, let $m_i$ have weight $\la_i$. Then $\chi_{-1}(V_i)$ is
obtained from the ordinary character $\chi(V_{\la_i})$ of the
irreducible $\g$--module $V_\la$, by replacing each $y_i^{\pm 1}$ with
$\bar{Y}_{i,a\psi(i)}^{\pm 1}$.
\end{lem}

\begin{proof}
According to \cite{L}, Proposition 33.2.3, the algebras
$\dot{U}^{\on{res}}_{-1} \G$ and $\dot{U}_{1}^{\on{res}} \G$ are
isomorphic, if $\g \neq \sw_{2n+1}$. Therefore the category of type 1
finite-dimensional $\dot{U}^{\on{res}}_{-1} \G$--modules is isomorphic
to the category of weight $\G$--modules. The statement of the lemma is
obtained by combining this fact with \thmref{spec char} and the
following result about $\chi_q$--characters, which follows from
\cite{FM}, Corollary 6.4:

Set $\varphi(i)=(1+\psi(i))/2$.
Suppose that $V$ is an irreducible $\uqg$--module, such that the
highest weight monomial $m$ in $\chi_q(V)$ belongs to
$\Z[Y_{i,a_0q^{\varphi(i)+2k}}]_{i\in I}^{k\in\Z}$.
Then all other
monomials in $\chi_q(V)$ belong to
$\Z[Y_{i,a_0q^{\varphi(i)+2k}}^{\pm1}]_{i\in I}^{k\in\Z}$.

In the case of $\g=\sw_{2n+1}$, we use the inclusion
$U_{-1}\widehat{\sw}_{2n+1}\to U_{-1}\widehat{\sw}_{2n+2}$ and Lemma
\ref{minus 1 lem} follows from Lemma \ref{comm diag for restr}.
\end{proof}

Combining this lemma with \thmref{hw of frob}, we obtain a complete
description of the $\ep$--characters of irreducible Frobenius
pull-backs in the case when $\ep^*=-1$.

These results together with Theorem \ref{decomp thm} reduce the
problem of calculating the $\ep$--characters of irreducible
finite-dimensional representations of $\ueg$ to the problem of
calculating $\ep$--characters of those representations which remain
irreducible when restricted to the subalgebra $U^{\on{fin}}_\ep \G$.

In particular, we obtain explicit formulas for the $\ep$--characters
of all irreducible finite-dimensional
$U^{\on{res}}_\ep\widehat\sw_2$--modules (cf. \cite{CP:root}, Theorem
9.6). We conjecture that in the case when $s>2$ the $\ep$--characters
of irreducible modules can be constructed by applying the algorithm
described in Section 5.5 of \cite{FM}.

\end{document}